\documentclass[letterpaper,10pt,conference]{IEEEtran}
\usepackage{amsmath,graphicx,amssymb}
\usepackage{multirow}
\usepackage{amsmath}
\usepackage[font=scriptsize,labelfont=bf]{caption}

\newtheorem{theorem}{Theorem}

\newtheorem{lemma}{Lemma}
\newtheorem{remark}{Remark}

\usepackage{epsfig}
\usepackage{subfigure}
\usepackage{geometry}
\title{Approximation learning methods of Harmonic Mappings in relation to Hardy Spaces }

\author{\IEEEauthorblockN{Zhulin Liu}
\IEEEauthorblockA{University of Macau\\
Faculty of Science and Technology\\
Taipa, Macao, China}

\and
\IEEEauthorblockN{C. L. Philip Chen}
\IEEEauthorblockA{University of Macau\\
Faculty of Science and Technology\\
Taipa, Macao, China}}

\begin{document}

\newgeometry{top=1in,bottom=0.75in,left=0.75in,right=0.75in}
%

 \maketitle
\begin{abstract}

A new Hardy space Hardy space approach of Dirichlet type problem based on Tikhonov regularization and Reproducing Hilbert kernel space is discussed in this paper, which turns out to be a typical extremal problem located on the upper upper-high complex plane. If considering this in the Hardy space, the optimization operator of this problem will be highly simplified and an efficient algorithm is possible. This is mainly realized by the help of reproducing properties of the functions in the Hardy space of  upper-high complex plane, and the detail algorithm is proposed. Moreover, harmonic mappings, which is a significant geometric transformation, are commonly used in many applications such as image processing, since it describes the energy minimization mappings between individual manifolds. Particularly, when we focus on the planer mappings between two Euclid planer regions, the harmonic mappings are exist and unique, which is guaranteed solidly by the existence of harmonic function. This property is attractive and simulation results are shown in this paper to ensure the capability of applications such as planer shape distortion and surface registration.  
\end{abstract}
\section{Introduction}
\label{sec:intro}

Geometric transformation has been widely used in many applications, such as image processing and pattern recognition. Among them, the harmonic geometric transformation plays a very important role. However, in the data acquisition process of various applications, shape distortion is common, which usually could be modeled by harmonic mappings.

A real-valued function mapping $u(x,y)$ in $R^2$ is harmonic if it satisfies Laplace's equation:
$$\Delta u = \frac{\partial^{2}u}{\partial x ^{2}} + \frac{\partial^{2} v}{\partial x ^{2}} =0.$$
It is important in many field of natural science, for the harmonic functions can describe the energy minimization mappings between individual manifolds. Recall that the Dirichlet problem for Laplace's equation consists of finding a solution $u$ on some domain $\Omega$ such that $u$ on the boundary of $\Omega$ is equal to some given function $g \in L^2(\partial \Omega)$. In this paper, a new Hardy space Hardy space approach of Dirichlet type problem based on Tikhonov regularization and Reproducing Hilbert kernel space is discussed. Unlike the other scholars who are interested in this classical problem, Saitoh(\cite{MaSa05}) has developed a practical approximation solution through the generalized inverses. The key point of his method is to approach the generalized inverse by a functional approximation in Sobolev reproducing Hilbert kernel spaces. Though the solution has solid mathematical background and graceful linear combination format, its accuracy rate is  highly related with the capacity of computers. It is possible to overcome this if we considering the original problem in another more convenient Hilbert space. 

It is necessary for us to introduce the definition of complex harmonic mappings, which relates closely with the planer harmonic mapping. A one-to-one mapping $u=u(x,y), v=v(x,y)$ from a region $\Omega$ in the $xy$-plane to a region $\Gamma$ in the $uv$-plane is a harmonic mapping if the two coordinate functions are harmonic (see \cite{CTM156}).  It is convenient to use the complex notation  $z=x+\imath y, w=u+ \imath v$, and to write $w=f(z)=u(z)+iv(z).$ Thus a complex-valued harmonic function is a harmonic mapping of a domain $\Omega \subset \mathbb{C}$ if and only if it is univalent in $\Omega$. Here $\mathbb{C}$ denotes the complex plane. The development of a theory of harmonic mappings in the plane has started as early as the 1920s, and more recently, complex analysts have been more interested in harmonic mappings as generalizations of conformal mappings.

In this research, based on complex analysis, we obtain a modified solution of Dirichlet problem and  try to construct  improved  models of shape distortion and surface registration combining with harmonic functions and mappings.

\section{Hardy Space Approach of Dirichlet Problem}
Our research begins with a modified solution of Dirichlet problem.
Without loss of generality, we only consider the situation that set $\Omega$ is located on the upper-high planer of $R^2$ which is equivalent to upper-high complex plane $\mathbb{C}^+$. Correspondingly, it is well known that the functions in $f(x+y\imath) \in H^2(\mathbb{C}^+)$ which can be written as below has the following properties.
  $$f(x+y\imath)=g(x,y)+h(x,y)\imath,$$ where $g(x,y), h(x,y)$ is harmonic functions.
Hence, it is natural for us modify and simplify a Dirichlet problem as follows: look for a function $f(x+y\imath) \in H^2(\mathbb{C}^+)$ such that
\begin{equation}
Re{f(x+y\imath)\rvert _{\partial \Omega}} =g(x,y)    \ \ \   on \ \ \ \partial \Omega.
\end{equation}
Here $\Re f$ is the real part the analytic function $f$. This is reasonable since in a simply connected domain $\Omega \subset \mathbb{C}$, a real-valued harmonic function $f$ has the representation $f=h+\overline{h}=Re{2h}$, where $2h$ is the analytic completion of $f$ , unique up
to an additive imaginary constant.
Hence, define an operator L :
\begin{equation*}
L: f(z)\in H^2(C^+) \mapsto \Re f\rvert_{\partial \Omega} \in L^2(\partial \Omega),
\end{equation*}
Therefore the Dirichlet problem can be rewritten to a linear equation of the form $$Lf=g.$$
It is obviously that for any $g\in L^2(\partial \Omega)$, it may not exist an exact function $f\in H^2(\mathbb{C}^+)$ that $Lf=g$. However, its generalized solution $f^*$ is available since the range of $L$ is dense in $L^2(\partial \Omega)$.
Finally it comes down to the following best approximation problem: Find an $f$ in $H^2(\mathbb{C}^+)$ which satisfies

\begin{align*}
&\text{argmin}_{f\in H^2(\mathbb{C}^+)}\|Lf-g\|^2 =\text{argmin}_{f\in H^2(\mathbb{C}^+)} \lVert Re f-g\rVert^2_{L^2(\partial \Omega)}.
\end{align*}

To make sure the existence of the extremal functions, we apply the Tikhonov regularization, the extremal problem can be considered as follows

\begin{equation}
\inf\limits_{f\in H^2(\mathbb{C}^+)}\left\{ \lambda \lVert f \rVert^2_{H^2(\mathbb{C}^+)} + \lVert Re f-g\rVert^2_{L^2(\partial \Omega)}\right\}.
\end{equation}
Here, $\lambda$ is the regularization parameters which usually fixes as a constant bigger than $0$.

Denote  Szeg\"{o} Kernel $K(z,w)$ in  $H^2(\mathbb{C}^+)$ as follows,
$$K(z,w)=\frac{\imath}{2 \pi }\frac{1}{z-\bar{w}}, \text{for}\ \ z,w \in \mathbb{C}^+.$$

To give a general conclusion of this kind of extremal problems, in the rest of this section, we simplify the formula by substituting the notation $Re$ by $L$ as below.
 $$\inf\limits_{F\in H^2(\mathbb{C}^+)}\left\{ \lambda \lVert F \rVert^2_{H^2(\mathbb{C}^+)} + \lVert LF-g\rVert^2_{L^2(\partial \Omega)}\right\}.$$

This regularization operator, however, can be reviewed as a modification of the norm of the reproducing Hilbert kernel space. It is important for us to describe the transferred spaces by deducing the new reproducing kernels, while the reproducing properties of the space could be ensured by the theories of reproducing Hilbert kernel space. The following theorem is a key part of the whole theory which can be checked in \cite{saitoh2005}. 

\begin{theorem}
In the reproducing Hilbert kernel space $H_K $ whose kernel is $K$, denote L: $H_K \mapsto \mathcal{H}$ as a
bounded linear operator, and an inner product is defined as 
$$\langle f_1, f_2\rangle_{H_{K_{\lambda}}}=\lambda \langle f_1, f_2\rangle_{H_K}+\langle Lf_1, Lf_2\rangle_{\mathcal{H}}$$
for $f_1,f_2 \in H_K$.  As a consequence, $(H_K,\langle \cdot \rangle_{H_{K_{\lambda}}}
)$ is also a reproducing kernel Hilbert space(RKHS). The new reproducing kernel is formulated by 
\begin{equation}
K_{\lambda}(p,q)=((\lambda + L^* L)^{-1} K_q)(p).
\end{equation}

Therefore, the following mapping which is constructed by the above inner production  
$$f\in H_K \mapsto \lambda \lVert f \rVert_{H_K} + \lVert Lf-g\rVert_{\mathcal{H}} $$
could reach the minimum and the minimum is attained only at $f_{g,\lambda} \in H_K$ such that
 \begin{equation}
(f_{g,\lambda})(p)= \langle g, LK_{\lambda}(\cdot,
p)\rangle_{\mathcal{H}}.
\end{equation}
Furthermore, if L is compact, the Moore-Penrose generalized inverse $f^*_{g}$ of L has an approximation expression, and the optimal solution of the extremal solution converges uniformly as following
$$\lim\limits_{\lambda \to 0} f_{g,\lambda}(p)=f^*_{g},$$ where $E_0$ is any subset of $E$ satisfying $\sup\limits_{p\in {E_0}}K(p,p)<\infty$.
\end{theorem}

Consequently, if substitute the notation $L$ by the operator $Re$ which is involved in Dirichlet problem, the solution of the extremal function can be represented below, 

\begin{align}
F^*_{\lambda,g}(p)&=\langle g, LK_{\lambda}(\cdot, p)\rangle_{L^2(\partial \Omega)} \notag\\
& =\langle g, Re{K_{\lambda}(\cdot, p)}\rangle_{L^2(\partial \Omega)} \notag\\
&=(\lambda + L^* L)^{-1}L^*g(p).
\end{align}

It is interesting that this formula is exactly same with the existing Tikhonov regularization solution which is the common method of dealing with ill posed problems \cite{tiki2}. This solution is proved to be convergent, and an error estimation $O(\lambda)$ is guaranteed in  \cite{tiki1}. As a result, we have a controlled convergence rates for our method.
However, for a general linear bounded operator $L$, it is not easy to achieve its conjugation operator effectively. Luckily, the operator $Re$ we discuss in this paper, has a friendly expression especially for the discrete data. The detail will be presented in the following section.


\section{Discrete Algorithm}

For practical applications, the data set we deal with is usually discrete. In this section we will give a description of discrete algorithm. Now for any given given real valued data $\{G(z_j)=G(x_j,y_j)=A_j\}, \text{where}  \ j=1,\dots,N$,  and $z_j=x_j+y_j\imath$ is the corresponding points of the boundary $\partial \Omega$ of $\mathbb{C}^+$. Here we treat the discrete data $g(z_j)$ as the real projection of some analytic function in $H^2(\mathbb{C}^+)$.

The following discretization is a rational $2$-norm discretization if the constant $\lambda$ is fixed,
\begin{equation}
\inf\limits_{F\in H^2(\mathbb{C}^+)}\{ \lambda \lVert F \rVert^2_{H^2(\mathbb{C}^+)} + \sum_{j=1}^N \lambda_j \lvert \text{Re}F(z_j)-A_j\rvert^2 \},
\end{equation}
where $\{\lambda_j\}_{j=1}^N (\lambda_j >0)$ are the weights parameters determined by the discretization.
As a result, the discrete norm of the new reproducing kernel space is reduced to the following norm square

\begin{equation}
\lambda \lVert F \rVert^2_{H^2(\mathbb{C}^+)} + \sum_{j=1}^N \lambda_j \lvert \text{Re}F(z_j)\rvert^2,
\end{equation}

Instead of approach the conjugation operator of the above norm, which is a normal process for other similar solutions, we will try to attain the new reproducing kernel defined by the above norm. It is possible by the following iteration steps.

We will start our algorithm by setting the boundary containing only $1$ point, i.e. $\Omega=\{z_1=x_1+y_1\imath\}$, the extremal problem should be modified as $$ \lambda \lVert F \rVert_{H^2(\mathbb{C}^+)} + \lambda_1 \lvert \text{Re}F(z_1)\rvert^2.$$ Therefore, the operator "L" is 
$$L: F(x+y\imath)\in H^2(\mathbb{C}^+) \mapsto \text{Re}F(z_1)=\text{Re}F(x_1+y_i\imath),$$
and the associated reproducing kernel could be deduced by the Fredholm integral equation of the second kind for concrete cases. The details of the deduction are omitted, while the solution expression is given by 

\begin{align}
K^{(1)}_{\lambda}(z,w)&=\frac{1}{\lambda} \cdot \notag \\
&(K(z,w)-\frac{\lambda_1 \text{Re}K(z_1,w)\text{Re}K(z,z_1)}{\lambda_1\text{Re}K(z_1,z_1)+\lambda}).
\end{align}
Where $K(z,w)$ is Szeg\"{o} kernel.

Next, denote the reproducing kernel of "L" determined spaces$H^2_{\lambda_1}(\mathbb{C}^+)$ as $K^{(1)}_{\lambda}(z,w)$. After putting a new point $z_2=x_2+y_2\imath$ to the boundary points, the new norm should be modified as 
 \begin{align*}
&\lambda \lVert F \rVert_{H^2(\mathbb{C}^+)} + \sum_{j=1}^2 \lambda_j \lvert \text{Re}F(z_j)\rvert^2  \\
&  = \lVert F \rVert_{H^2_{\lambda_1}(\mathbb{C}^+)}+ \lambda_2 \lvert \text{Re}F(z_2)\rvert^2.
 \end{align*}
Continue solving the Fredholm integral equation of the second kind we would get the step $2$ kernel $K^{(2)}_{\lambda}(z,w)$ associated with the updating operator $L$. 
\begin{align}
K^{(2)}_{\lambda}(z,w)&==K^{(1)}(z,w)\notag\\
&-\frac{\lambda_2 \text{Re}K^{(1)}(z_2,w)\text{Re}K^{(1)}(z,z_2)}{\lambda_2\text{Re}K^{(1)}(z_2,z_2)+1}.
\end{align}

After all the $N$ points are added into the boundary set $\Omega$, we have that 

\begin{align}
K^{(n)}_{\lambda}(z,w)&=K^{(n-1)}(z,w)\\
&-\frac{\lambda_2  \text{Re} K^{(n-1)}(z_{n},w) \text{Re} K^{(n-1)}(z,z_{n-1})}{\lambda_n \text{Re} K^{(n-1)}(z_{n-1},z_{n-1})+1},
\end{align}
Here $K^{(N)}_\lambda$ is the reproducing kernel we need for the given extremal problems.
As a result, the approximation of the solution function could be attained by the inner production with $\{G(z_j)=G(x_j,y_j)=A_j\}, \text{where}  \ j=1,\dots,N$ 
\begin{align}
 F_{\{A_j, z_j\},\lambda}(x,y)&==\sum_{j=1}^N A_j \lambda_j\text{Re}K^{(N)}_\lambda(x+y\imath,x_j+y_j\imath).
\end{align}


Now we have an analytic expression of the extremal function involved with the iteration kernel $K^{(N)}_\lambda$. It is more connivent if we could deduce an clearer expression presented with the given data. Therefore we try to rewrite the algorithm in matrix form.

Denote
$$\hat{K}=(K(z,z_1),\cdots,K(z,z_N))^T$$

$$ \hat{K}^{(n)}_{\lambda}=(K^{(n)}_{\lambda}(z,z_1),\cdots,K^{(n)}_{\lambda}(z,z_N))^T.$$
and denote
\begin{equation}
A^{(1)}=\frac{1}{\lambda }I-\left [ a^T|0 \right ],
\end{equation}
where
$$a=\begin{pmatrix}\frac{\lambda_1\text{Re} K(z_1,z_1)}{\lambda(\lambda_1\text{Re} K(z_1,z_1)+\lambda)}&\cdots&\frac{\lambda_1\text{Re} K(z_1,z_N)}{\lambda(\lambda_1\text{Re} K(z_1,z_1)+\lambda)}\end{pmatrix}$$
We can easily get that 
\begin{equation}
\text{Re} \hat{K}^{(1)}_{\lambda}=A^{(1)} \text{Re} \hat{K}.
\end{equation}
Similarly, we have
\begin{equation}
\text{Re}\hat{K}^{(n)}_{\lambda}=A^{(n)} \text{Re} \hat{K}^{(n-1)}_{\lambda}.
\end{equation}

Finally, for the $N$th iteration result $\text{Re}\hat{K}^{(N)}_{\lambda}$, we have that
\begin{align}
\text{Re}\hat{K}^{(N)}_{\lambda} &=A^{(N)} \text{Re} \hat{K}^{(N-1)}_{\lambda}=\prod\limits_{n=1}^{N}A^{(n)}\text{Re} \hat{K}\notag\\
&=\hat{A}\text{Re} \hat{K}.
\end{align}
Here $\hat{A}=\prod\limits_{n=1}^{N}A^{(n)}.$

For the extremal solution functions, by substituting $\text{Re} K^{(N)}_{\lambda}$ in $(13)$ by the above formula, we have that  

\begin{align}
F_{\{A_j, z_j\},\lambda}(z)&= A'\hat{A}\text{Re} \hat{K}\notag\\
& =\sum_{j=1}^N c_j \text{Re} K(z,z_j).
\end{align}
Where $A'=(\lambda_1 A_1,\lambda_2 A_2,\cdots,\lambda_N A_N)$. The parameters $c_j$ are easily to calculate since they are combinations of the parameters $\{\lambda_i\}, \lambda,$ and $A_j$.

Hence for the given arbitrary boundary points $\Omega$, we conclude that the approximation of Dirichlet problems could actually be viewed as the linear combinations of the discrete Szeg\"{o} Kernels. It is potential for the practical applications. Now we have the following theorem.
%
%
\begin{theorem}
Denote $\{G(z_j)=G(x_j,y_j)=A_j\}, \text{where}  \
j=1,\dots,N$,  and $z_j=x_j+y_j\imath$ as a sequence of real valuled points on the
boundary $\partial \Omega$ of $C^+$,  and $\lambda, \{\lambda_i\} >0, j= 1,\dots,N$, $F_{\{A_j,
z_j\},\lambda}(z)$ are artriburary fixed values, then the function $F_{\{A_j,z_j\},\lambda}(z)$  in $(18)$ is an approach of the given discrete Dirichlet problem in the upper high complex Hardy space.
Moreover, the Moore-Penrose generalized inverse could also be achieved by the following limitation.
$$\lim\limits_{\lambda \to 0} F_{\{A_j, z_j\},\lambda}(z)=F^*_{\{A_j\}},$$
which means $F^*_{\{A_j\}}$ is the solution of the problem
$$\Delta u=0, u(x_j,y_j)=A_j, \ \ j=1,\dots,N$$

\end{theorem}

\begin{remark}
The constraint that the boundary function of the Dirichlet problems should be the real part of some analytic function is not critical for the discrete part of our algorithm. For the approximation approach in Hardy space always exists for the discrete data.
\end{remark}
\begin{remark}
It is one noteworthy property that our algorithm is not limited to the regular domain, which is usually a disadvantage of other methods. Moreover, we do not need any carefully  divided grids, which is an attractive advantage especially in shape distortion.
\end{remark}

\section{Shape Distortion}
Generally, nonlinear linear shape distortions are more challenging topic in practical applications. In our research, we mainly deal with the shape distortion in image processing, and the model can be described as below.
To construct a harmonic models of nonlinear shape distortion, we have the following assumptions in \cite{tang, harmonic surfaces}. Consider a nonlinear transformation $T$ :
$$T:(\xi,\eta)\in \Omega \rightarrow (x,y)\in \Gamma, $$
where $\Omega$ is the image in $\xi O \eta$ with the boundary $\partial \Omega$, and $\Gamma$ is the image in $X O Y$ with the boundary $\partial \Gamma$.
To construct a stable mapping between $\Omega$ and $\Gamma$, it is reasonable to assume that the transformation is harmonic.On the basis of the above ideas, we develop the harmonic models as follows:
$$T:(\xi,\eta)\in \Omega \rightarrow (x,y)\in \Gamma,$$
where
$$\Delta x(\xi,\eta) = 0, x|_{\partial \Omega} = f(\xi,\eta),$$
$$\Delta y(\xi,\eta) = 0, y|_{\partial \Omega} = g(\xi,\eta).$$

Under this assumption, we have the following lemma.
\begin{lemma}
Let $$T:(\xi,\eta)\in \Omega \rightarrow (x,y)\in \Gamma, $$ be one-to-one, and $T \in C^1(\Omega), f\in C^0(\Omega).$ Then,

$$\int\int_{\Gamma_i} f(x,y) dx dy\\
= \int\int_{\Omega_i} f(x(\xi,\eta),y(\xi,\eta)) J(\xi,\eta)d\xi d\eta,$$
where $\Omega_i \in \Omega$, $\Omega_i\rightleftarrows\Gamma_i$ and the Jacobian determinant
$$J(\xi,\eta)=
\begin{vmatrix}
\frac{\partial x}{\partial \xi} & \frac{\partial y}{\partial \xi}\\
\frac{\partial x}{\partial \eta} &\frac{\partial y}{\partial \eta}
\end{vmatrix}$$
\end{lemma}
On the basis of the above lemma, we will develop the splitting-integral methods for $T$, the splitting-shooting methods for $T^{-1}$, and their theoretical analysis.
Firstly by adopting the trivial technology in (\cite{Tintegrating},\cite{Tshooting}), we have a simple example below.
\begin{figure}
  \centering
  \begin{minipage}[t]{.3\linewidth}
 \includegraphics[width=3cm]{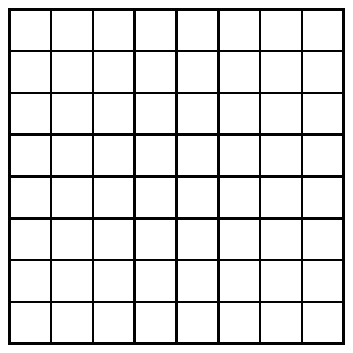}
  \caption{The original $8\times8$ square grid}
  \end{minipage}
  \begin{minipage}[t]{.3\linewidth}
  \includegraphics[width=3cm]{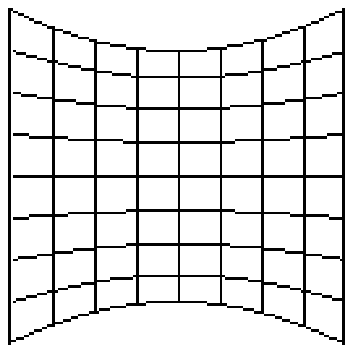}
  \caption{Shape distorton under a quadratic pressing}
  \end{minipage}
  \begin{minipage}[t]{.3\linewidth}
 \includegraphics[width=3cm]{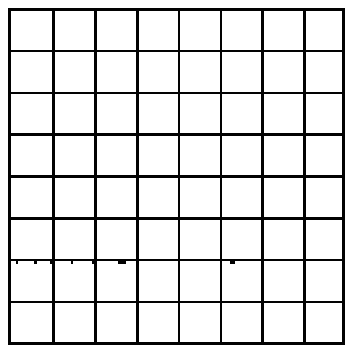}
  \caption{Recovered image from the distortion image}
  \end{minipage}
\end{figure}

\section{Simulation results}

Consider the girl-image  that consists of $256\times256$ pixels of $256$ greyness levels, and choose the piecewise constant interpolation we have mentioned before. The results show that this model is valid for the harmonic geometry transformations. It is clear that the image is recovered corrected and the details of the image are remained.

\begin{figure}
  \centering
  \begin{minipage}[t]{.35\linewidth}
 \includegraphics[width=2cm]{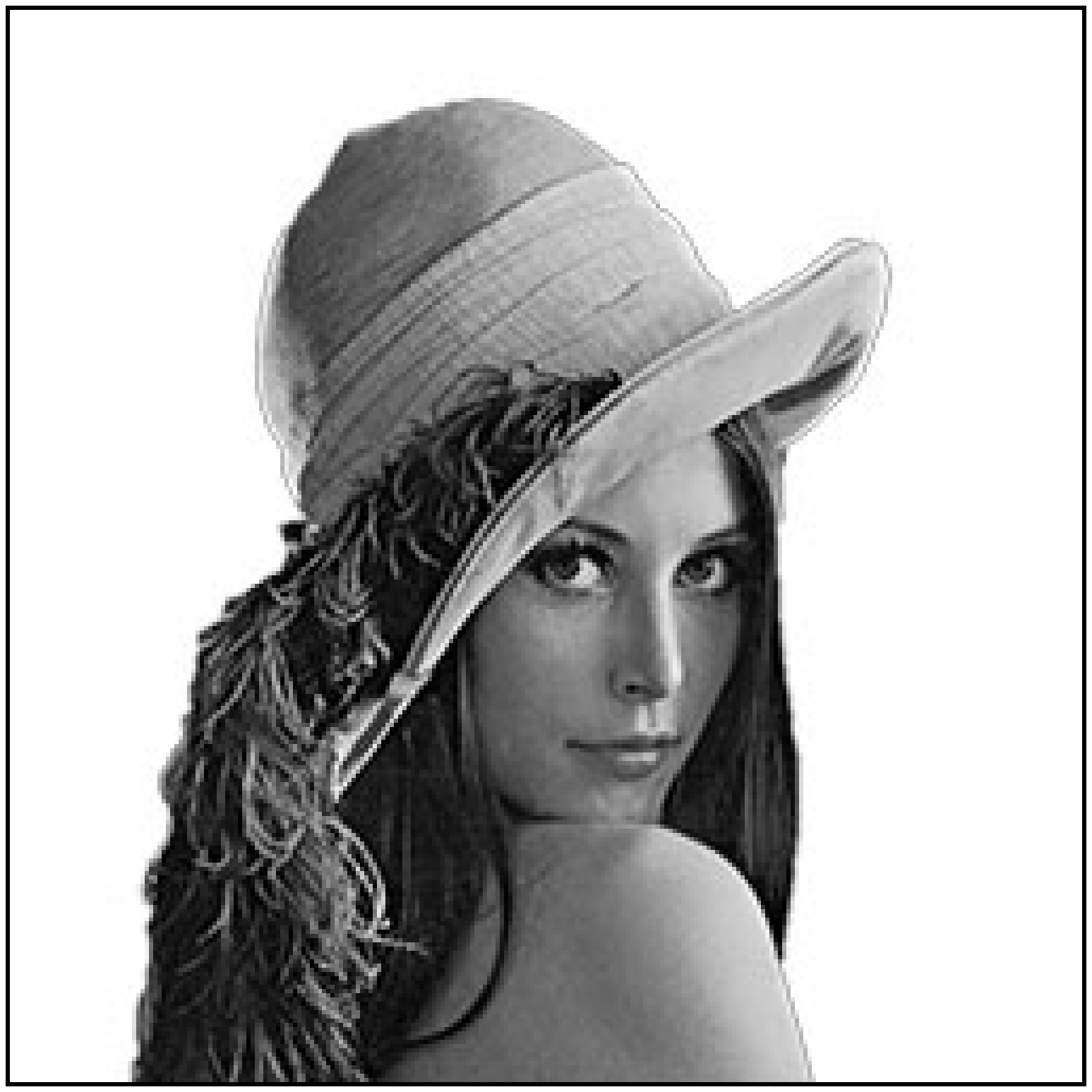}
  \caption{Girl original image of $256\times256$ with $256$ greyness levels}
  \end{minipage}
  \begin{minipage}[t]{.35\linewidth}
  \includegraphics[width=2cm]{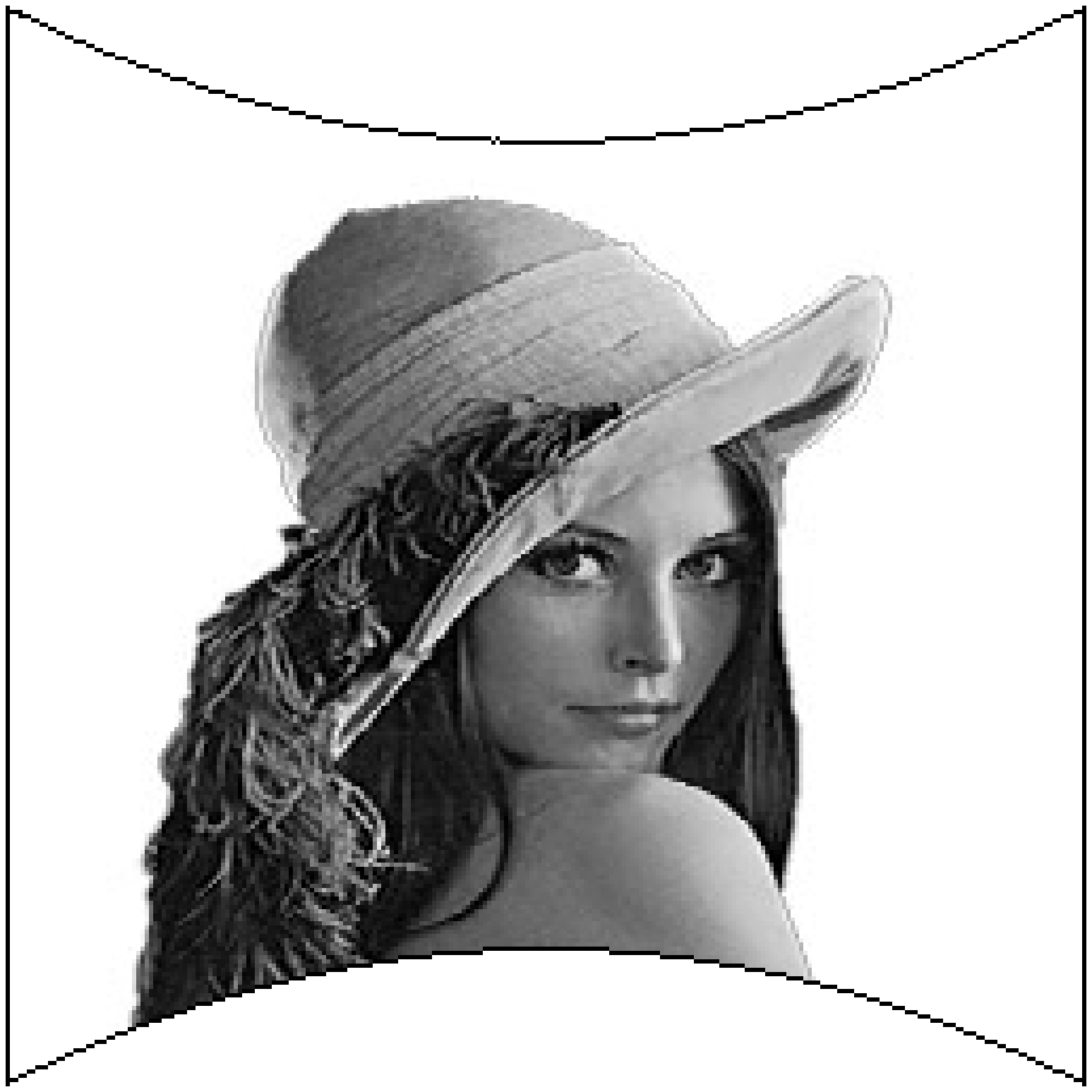}
  \caption{Distorted girl image of $256\times256$ under $T$}
  \end{minipage}
  \begin{minipage}[t]{.35\linewidth}
 \includegraphics[width=2cm]{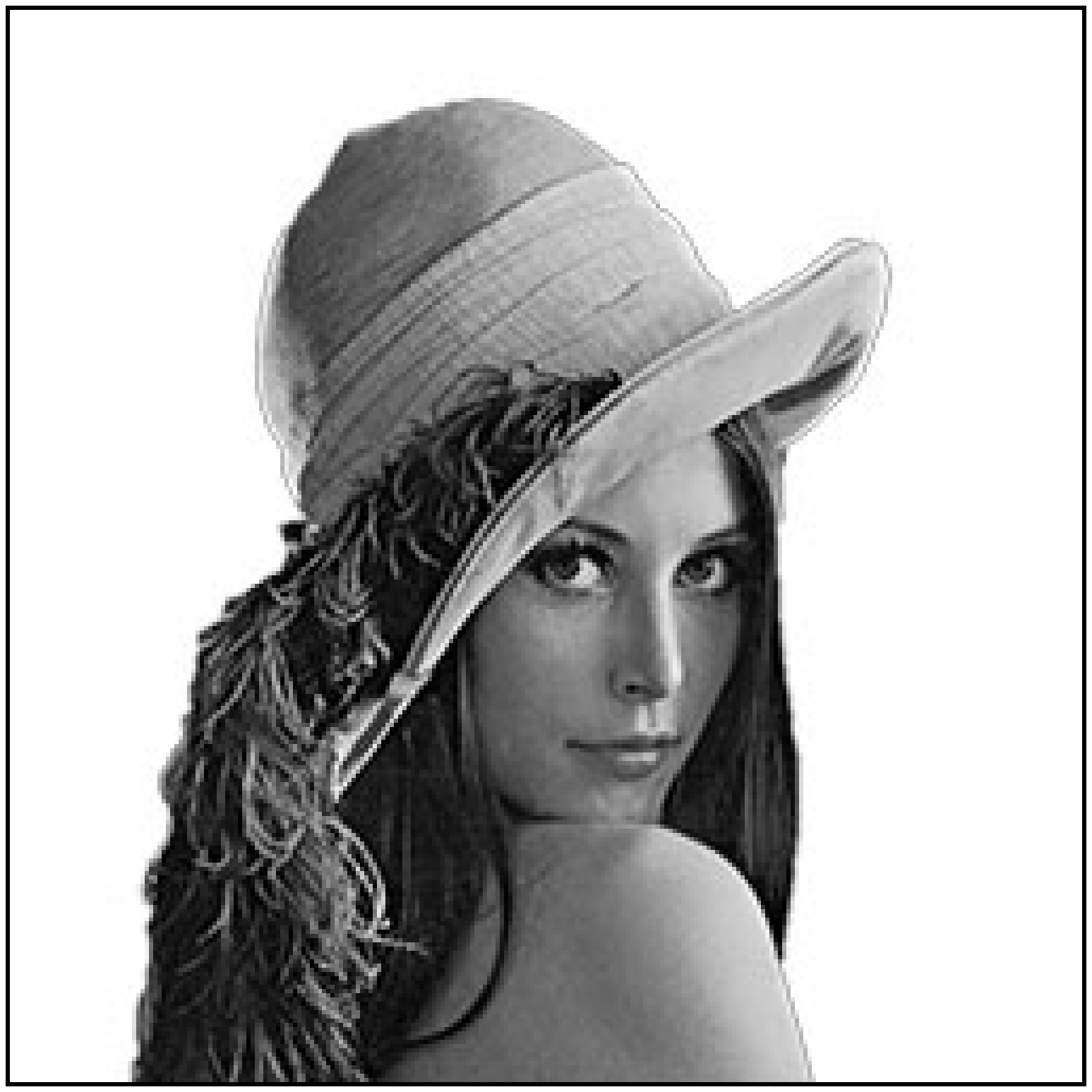}
  \caption{{Girl image under $T^{-1}T$ by discrete methods}}
  \end{minipage}
\end{figure}

\section{Conclusion}

 in this paper, a new Hardy space Hardy space approach of Dirichlet type problem based on Tikhonov regularization and Reproducing Hilbert kernel space is discussed, which turns out to be a typical extremal problem located on the upper upper-high complex plane.  This is mainly realized by the help of reproducing properties of the functions in the Hardy space of  upper-high complex plane, and the detail algorithm is proposed. Finally, it is found that under the assumption we give, the optimal solution is actually the linear combination of Szeg\"{o} kernels. Moreover,our algorithm is not limited to the regular domain and carefully divided domains. This is potential for some practical applications. 
  
 Harmonic mappings, which is a significant geometric transformation, are commonly used in many applications such as image processing, since it describes the energy minimization mappings between individual manifolds. Particularly, for planer shape distortion, since the planer mappings between two Euclid planer regions are exist and unique, which is guaranteed solidly by the existence of harmonic function, an algorithm scheme is proposed for most kind of harmonic based shape distortion. This property is attractive and simulation results are shown in this paper to ensure the capability of applications. 

\section{acknowledge}

This research is supported, in part, by Macao Science and Technology Development Fund under Grant No.: 019/2015/A1, and Research Committee of the University of Macau under Grant No.: RDG008/FST-CCL/2012, and the National Nature Science Foundation of China under Grant No.: 61572540.


\end{document}